\documentclass[12pt]{article} 

\usepackage{enumerate}

\usepackage{amsmath}
\usepackage{amsfonts}
\usepackage{amssymb}
\usepackage{url}

\DeclareMathAlphabet{\eufrak}{U}{}{}{} 
\SetMathAlphabet\eufrak{normal}{U}{euf}{m}{n}
\SetMathAlphabet\eufrak{bold}{U}{euf}{b}{n}
\numberwithin{equation}{section}

\newenvironment{Proof}{\removelastskip\par\medskip
\noindent{\em Proof.} \rm}{\penalty-20\null\hfill$\square$\par\medbreak}

\newenvironment{Proofy}{\removelastskip\par\medskip
\noindent{\em Proof } \rm}{\penalty-20\null\hfill$\square$\par\medbreak}

\allowdisplaybreaks

 \def\G{{\mathord{{\rm {\sf G}}}}}
 \def\g{{\mathord{{\rm {\sf g}}}}}

 \def\inte{{\mathord{\mathbb N}}}
 
 \def\qu{{\mathord{\mathbb Z}}}

 \def\Var{{\mathrm{{\rm Var}}}}

 \def\div{{\mathrm{{\rm div}}}}

 \def\inte{{\mathord{{\rm I\kern-3pt N}}}}
 \def\sZZ{{\rm Z\kern-.45em{}Z}}
 \def\sQQ{{\kern 0.27em \vrule height1.45ex width0.03em depth0em
           \kern-0.30em \rm Q}}
 \def\qu{{\mathchoice
         {\sQQ}
         {\sQQ}
   {\kern 0.225em \vrule height1.05ex width0.025em depth0em \kern-0.25em \rm Q}
   {\kern 0.180em \vrule height0.78ex width0.020em depth0em \kern-0.20em \rm Q}
         }}
 \def\sGG{{\kern 0.27em \vrule height1.45ex width0.03em depth0em
           \kern-0.30em \rm G}}
 \def\gg{{\mathchoice
         {\sGG}
         {\sGG}
   {\kern 0.225em \vrule height1.05ex width0.025em depth0em \kern-0.25em \rm G}
   {\kern 0.180em \vrule height0.78ex width0.020em depth0em \kern-0.20em \rm G}
         }}

 \newtheorem{prop}{Proposition}[section]
 \newtheorem{lemma}[prop]{Lemma}
 \newtheorem{definition}[prop]{Definition}
 \newtheorem{corollary}[prop]{Corollary}

 \def\div{{\mathrm{{\rm div}}}}

\def\E{\mathop{\hbox{\rm I\kern-0.20em E}}\nolimits}

\newcommand{\dee}{\mbox{$I  \! \! \! \, D$}}
\def\dee{{\mathord{\mathbb D}}}

\newcommand{\Supp}{\mathrm{\tiny Supp~}}

\def\inte{{\mathord{\mathbb N}}}

\def\qu{{\mathord{\mathbb Q}}}

\def\PP{{\mathord{{\rm I\kern-2.8pt P}}}}

\def\E{\mathop{\hbox{\rm I\kern-0.20em E}}\nolimits}

\def\Var{\mathop{\hbox{\rm Var}}\nolimits}

\def\div{{\mathord{{\rm  \hskip0.1cm div }}}}

\def\inte{{\mathord{{\rm I\kern-2.8pt N}}}}
\def\PP{{\mathord{{\rm I\kern-2.8pt P}}}}

\def\real{{\mathord{\mathbb R}}}

\def\inte{{\mathord{\mathbb N}}}

\def\qu{{\mathord{\mathbb Z}}}

\def\Var{{\mathrm{{\rm Var}}}}

\usepackage[colorlinks=true, urlcolor=blue,linkcolor=blue, citecolor=blue]{hyperref}

 \newcounter{hyp}
\title{\huge 
 Stein approximation for multidimensional Poisson random measures by third cumulant expansions 
 } 
 
\author{
\large 
Nicolas Privault 
\\ 
\normalsize 
Division of Mathematical Sciences 
\\ 
\normalsize 
School of Physical and Mathematical Sciences 
\\ 
\normalsize 
Nanyang Technological University 
\\ 
\normalsize 
21 Nanyang Link 
\\ 
\normalsize 
Singapore 637371
\\ 
\normalsize 
nprivault@ntu.edu.sg
}

\usepackage{fixltx2e,amsmath}
\MakeRobust{\eqref}

\begin{document}
 
\maketitle 
 
\vspace{-0.67cm} 

\begin{abstract} 
We obtain Stein approximation bounds for stochastic integrals with respect to a Poisson random measure over $\real^d$, $d\geq 2$. This approach relies on third cumulant Edgeworth-type expansions based on derivation operators defined by the Malliavin calculus for Poisson random measures. The use of third cumulants can exhibit faster convergence rates than the standard Berry-Esseen rate for some sequences of Poisson stochastic integrals. 
\end{abstract} 

\baselineskip0.7cm 
 
\section{Introduction} 
Stein approximation bounds for stochastic integrals
with respect to a Poisson random measure have been obtained in
\cite{utzet2} using finite difference operators on the
Poisson space.
In this paper we derive related bounds for compensated Poisson stochastic integrals
$\delta ( u ) := \int_{\real^d} u_x ( \gamma ( dx ) - \lambda (dx))$
of compactly supported processes $(u_x)_{x\in \real^d}$ with respect to a Poisson
random measure $\gamma (dx)$ with intensity 
the Lebesgue measure $\lambda (dx)$ on $\real^d$, $d\geq 2$.
In contrast with \cite{utzet2},
our approach is based on derivation operators
and Edgeworth type expansions that involve the
third cumulant of Poisson stochastic integrals,
and can result into faster convergence rates,
see e.g. \eqref{dlf} below. 
\\
 
 Edgeworth type expansions have been obtained on the Wiener space in 
 \cite{nourdinpeccati2009}, 
 \cite{campese}, by a construction of cumulant operators 
 based on the inverse $L^{-1}$ of the
 Ornstein-Uhlenbeck operator, 
 extending the results of \cite{nourdinpeccati} 
 on Stein approximation and Berry-Esseen bounds. 
\\ 
 
In Proposition~\ref{djklsad1} we derive Edgeworth type expansions of the form 
\begin{equation}
\label{otf} 
 E\left[ 
 \delta ( u ) 
 {{g}} ( \delta ( u ) ) 
 \right] 
 = 
 E\big[ 
 \Vert u \Vert_{L^2(\real^d)}^2 {{g}}' ( \delta ( u ) ) 
 \big] 
 +
 \sum_{k=2}^n 
 E\left[ 
 {{g}}^{(k)} ( \delta ( u ) ) 
 \Gamma_{k+1}^u {\bf 1} 
 \right]
 +
 E\left[ 
 {{g}}^{(n+1)} ( \delta ( u ) ) 
 R_n^u
 \right]
\end{equation} 
 when the random field $(u_x)_{x\in \real^d}$ is predictable with respect to 
 a given total order on $\real^d$, where $\Gamma_k^u$ is a cumulant
 type operator and $R_n^u$ is a remainder term,  
 defined using the derivation operators
 of the Malliavin calculus on the Poisson space.
\\ 

 Based on \eqref{otf}, in Corollary~\ref{djkd} we deduce Stein
 approximation bounds of the form
\begin{eqnarray} 
\nonumber 
 d ( \delta ( u ) , {\cal N} ) 
&  \leq &  
 | 1 - \Var [ \delta (u) ] | 
 + 
 \sqrt{ 
 \Var \big[ 
 \Vert u \Vert_{L^2(\real^d )}^2 
 \big] 
 }
 \\
 \nonumber
 & &
 + 
   E\left[ 
   \left|
   \int_{\real^d } u_x^3 \ \! \lambda ( dx ) 
 + 
\left< 
 u , 
 D 
 \int_{\real^d } u^2_x \ \! \lambda ( dx ) 
 \right>_{L^2 (\real^d )}   
 \right|
 \right] 
+ 
   E\left[ 
 | 
 R_1^u
 | 
 \right] 
, 
\end{eqnarray} 
 where $D$ is a gradient operator acting on Poisson functionals,
 and ${\cal N}\simeq {\cal N}(0,1)$ is a standard
 Gaussian random variable,
 see also Proposition~\ref{djhlkdsa}.
 Here, 
$$ 
 d(F,G) := 
 \sup_{h\in {\cal L} }|\mathrm{E}[h( F )]-\mathrm{E}[h( G )]| 
$$
 is the Wasserstein distance between
 the laws of two random variables $F$ and $G$, 
 where ${\cal L}$ denotes the class of $1$-Lipschitz functions
 on $\real$. 
\\

 In particular, when $f$ is a differentiable deterministic function with
 support in the closed centered ball $B(R):=B(0;R)$ with radius $R>0$
 we obtain bounds of the form 
\begin{eqnarray} 
\label{jkdf} 
 d \left( \int_{\real^d} f(x) ( \gamma ( dx ) - \lambda (dx) ) , {\cal N} \right) 
 & \leq & 
 \big|
 1
 -
 \Vert
 f
 \Vert_{L^2(\real^d )}^2
 \big|
 +
 \left|
 \int_{\real^d} f^3(x) \ \! \lambda ( dx ) 
 \right|
 \\
 \nonumber
 & &
 + 
 8 ( K_d v_d R )^2 
  \Vert f \Vert_{L^2 (\real^d )}
  \Vert \nabla^{\real^d} f \Vert_{L^\infty (\real^d;\real^d)}^2,
\end{eqnarray} 
where $v_d$ denotes the volume of the unit ball in $\real^d$
and $K_d>0$ is a constant depending only on $d\geq 2$. 
 The bound \eqref{jkdf} can be compared to the classical Stein bound
\begin{equation} 
 \label{e1}
 d \left( \int_{\real^d} f(x) ( \gamma ( dx ) - \lambda (dx) ) , {\cal N} \right) 
 \leq
 \big|
 1
 -
 \Vert
 f
 \Vert_{L^2(\real^d )}^2
 \big|
 +
 \int_{\real^d} | f^3(x)| \ \! \lambda ( dx ),
 \end{equation} 
for compensated Poisson stochastic integrals,
see Corollary~3.4 of \cite{utzet2},
which involves the $L^3(\real^d)$ norm of $f$ instead
of third cumulant $\kappa_3^f = \int_{\real^d} f^3(x) \ \! \lambda ( dx ) $
of $\int_{\real^d} f(x) ( \gamma ( dx ) - \lambda (dx) )$,
and relies on the use of finite difference operators, 
see Theorem~3.1 of \cite{utzet2} and \S~4.2 of
\cite{bourguin}.
\\ 

For example when $f_k$, $k\geq 1$,
is a radial function given on $B(Rk^{1/d})$ by 
$$
f_k (x) := \frac{1}{C \sqrt{k} } g\left(\frac{|x|_{\real^d}}{k^{1/d}}\right)
, 
$$
where $g\in {\cal C}^1_0 ([0,R])$ is continuously differentiable on $[0,R]$,
and 
$$
C^2 := \int_0^R g^2(r) r^{d-1} dr < \infty,
$$
so that $\Vert f_k \Vert_{L^2(B(Rk^{1/d}))}=1$,
 the bound \eqref{e1} yields the standard Berry-Esseen
 convergence rate 
$$ 
 d \left( \int_{B(Rk^{1/d})} f_k(x) ( \gamma ( dx ) - \lambda (dx) ) , {\cal N} \right) 
 \leq 
\frac{v_d}{C^3 \sqrt{k}} \int_0^R |g (r)|^3 r^{d-1} dr,
\qquad k \geq 1. 
$$
 While \eqref{jkdf} does not improve on \eqref{e1} 
 when the function $f$ has constant sign, if 
 $g$ satisfies the condition 
$$
  \int_0^R g^3 (r) r^{d-1} dr =0, 
$$
then the third cumulant bound
\eqref{jkdf} yields the $O(1/k)$ convergence rate
\begin{equation}
  \label{dlf} 
 d \left( \int_{B(Rk^{1/d})} f_k(x) ( \gamma ( dx ) - \lambda (dx) ), {\cal N} \right) 
 \leq
 \frac{ 2 ( 2 K_d v_d R )^2 d }{ k C^2}\Vert g' \Vert_\infty^2, 
\end{equation} 
which improves on the standard Berry-Esseen rate, see Section~\ref{s4}. 
\\
 
In Sections~\ref{s2} and \ref{s2-1} we recall some background material
on the Malliavin calculus and differential geometry
on the Poisson space, by revisiting the approach
of \cite{pratprivault}, \cite{prconnection} using the recent
constructions of \cite{acosta} and references therein
on the solution of the divergence problem. 
In Section~\ref{s3} we derive Edgeworth type expansions
for the compensated Poisson stochastic integral
$\delta (u)$,
based on a family of cumulant operators that are associated to the
random field $(u_x)_{x\in \real^d}$.  
 In Section~\ref{s4} we derive Stein type approximation bounds for stochastic integrals,
  with deterministic examples. 
\\
  
While this paper is dealing with Poisson random measures
on $\real^d$ with $d\geq 2$,
the special case $d=1$ requires a different treatment 
for the standard Poisson process on the real half line,
see \cite{privaultthird}, and the $d$-dimensional setting of the present paper
shows significant differences with the one-dimensional case.
\subsection*{Preliminaries}
Let $d\geq 2$ and $0<R<R':=2R$.
We recall the existence of a ${\cal C}^\infty$ kernel function
$\G_\eta :B(R')\times B(R')\to \real^d$
defined as
$$
\G_\eta (x,y) : = \int_0^1 \frac{(x-y)}{s} \eta \left( y + \frac{x-y}{s} \right)
\frac{ds}{s^d},
\qquad
x,y \in B(R'),
$$
where $\eta\in {\cal C}^\infty_0 (B(R'))$ is such that 
$\int_{B(R)} \eta(x) dx = 1$, see \cite{acosta}, and satisfying the following properties: 
\begin{enumerate}[i)] 
\item The kernel $\G_\eta (x,y)$ satisfies the bound 
\begin{equation}
\label{gxy1} 
|\G_\eta (x,y)|_{\real^d} \leq \frac{K_d}{|x-y|^{d-1}_{\real^d}},
\qquad
x,y \in B(R'), 
\end{equation} 
for a constant $K_d>0$ depending only on $d$, 
see Lemma~2.1 of \cite{acosta},
by choosing $K_d$ and the function $\eta \in {\cal C}^\infty_c(B(R'))$ therein
so that $\Vert \eta\Vert_\infty \leq (d-1) K_d (R')^{-d}$. 
\item For any $p>1$ and 
  $g\in L^p(B(R'))$
  the function
\begin{equation}
\nonumber 
f(x) : = \int_{B(R')} \G_\eta (x,y) g(y) \ \! \lambda (dy), \qquad x\in B(R'),
\end{equation} 
 satisfies the bound 
\begin{equation}
\label{gxy3} 
\Vert f \Vert_{L^p(B(R');\real^d)}
\leq K_d v_d R' \Vert g \Vert_{L^p(B(R'))}, \qquad p>1, 
\end{equation} 
which follows from Young's inequality and \eqref{gxy1},
cf. Theorem~2.4 in \cite{acosta}. 
\item For any $h\in {\cal C}^\infty_0 (B(R'))$
  we have the relation 
\begin{equation}
\label{gxy4} 
h(y) - \int_{B(R')\setminus B(R)} h(x) \eta(x) \ \! \lambda (dx) =  \int_{B(R')} \langle \G_\eta (x,y) , \nabla^{\real^d}_x h(x) \rangle_{\real^d} \lambda (dx),
\quad
y\in B(R'), 
\end{equation} 
cf. Lemma~2.2 in \cite{acosta}, by taking
$\eta\in {\cal C}^\infty_c(B(R')\setminus B(R))$.
In particular, when $h\in {\cal C}^\infty_0 (B(R))$ we have
\begin{equation}
\label{gxy4-2} 
h(y) =  \int_{B(R')} \langle \G_\eta (x,y) , \nabla^{\real^d}_x h(x) \rangle_{\real^d}
\lambda (dx), \qquad y\in B(R'). 
\end{equation} 
\end{enumerate} 
\noindent 
 An extension of the framework of this paper
 by replacing $B(R)$ with a compact $d$-dimensional Riemannian manifold $M$ 
 and $\lambda (dx)$ with the volume element of $M$
 requires the Laplacian ${\cal L} = \div^M\nabla^M$ to be
 invertible on ${\cal C}_c^\infty (M)$, with  
$${\cal L}^{-1} u (x) = \int_M \g (x,y) u(y)\ \! \lambda (dy),
 \qquad x\in M, \ u\in {\cal C}^\infty_c (M), 
$$ 
 where $\g (x,y)$ is the heat kernel on $M$.
 In this case we can define $\G_\eta (x,y) \in \real^d$ as 
$$
 \G_\eta (x,y) = \nabla^M_x \g (x,y), \qquad {\lambda} \otimes {\lambda} (dx,dy)-a.e. 
$$ 
 with the relation
$$
 \nabla^M_x {\cal L}^{-1} u (x) 
 = \int_M u (y) \G_\eta (x,y) \ \! \lambda (dy)\in T_xM, \qquad x\in M, \ u \in {\cal C}^\infty_c (M),
$$ 
 from which the divergence inversion relation \eqref{gxy4-2} holds
 by duality. 
\section{Gradient, divergence and covariance derivative} 
\label{s2} 
There exists different notions of gradient and divergence
operators for functionals of Poisson random measures.
The operators of \cite{akr}, \cite{privaulttorrisi}, \cite{clausel},
and their associated integration 
by parts formula rely on an $\real^d$-valued gradient
for random functionals and a divergence operator which is
associated to the non-compensated Poisson stochastic 
integral of the divergence of $\real^d$-valued random fields.
This particularity, together with a lack of a suitable commutation relation
between gradient and divergence operators on Poisson functionals,
makes this framework difficult to use for a direct analysis of Poisson stochastic
integrals, while it has found applications to statistical estimation
and sensitivity analysis \cite{clausel}, \cite{privaulttorrisi}. 
\\
 
In this paper we use the construction of \cite{pratprivault},
\cite{prconnection} which relies on real-valued tangent processes
and of a divergence operator that directly
extends the compensated Poisson stochastic integral.
This framework also allows for simple commutation relations 
between gradient and divergence operators using
the deterministic inner product in $L^2(\real^d,\lambda)$,
see Proposition~\ref{commut}, and it naturally involves
the Poisson cumulants, see Definition~\ref{def01} and Relation~\eqref{ifinad2}. 
\subsubsection*{Gradient operator}  
In the sequel we consider a Poisson random measure
$\gamma (dx )$ on a probability space
$(\Omega , {\cal F} , P)$ and we 
let $\{X_1,\ldots ,X_n\}$ denote the configuration points
 of $\gamma (dx )$ when $B(R)$ contains $n$ points in the configuration $\gamma$, i.e.
 when $\gamma (B(R))=n$.
\begin{definition} 
 Given $A$ a closed subset of $B(R')$, we let ${\cal S}_A$ denote the
set of random functionals $F_A$ of the form 
\begin{equation}
  \label{f} 
 F_A =
  \sum_{n=0}^\infty 
  {\bf 1}_{\{ \gamma ( B(R))=n\} } 
 f_n \left( X_1,\ldots , X_n \right)
 ,
\end{equation} 
where
 $f_0\in \real$ and $(f_n)_{n\geq 1}$ is a sequence 
of functions 
satisfying the following conditions:
\begin{enumerate}[-]
\item
  for all $n \geq 1$, $f_n \in {\cal C}^\infty_c ( A^n )$
  is a symmetric function in $n$ variables, 
\item
 for all $n\geq 1$ and $i=1,\ldots , n$ we have the continuity condition
\begin{equation}
   \label{cont} 
 f_n \left( x_1,\ldots , x_n \right)
 =
 f_{n-1} \left( x_1,\ldots , x_{i-1},x_{i+1}, \ldots , x_n \right)
 ,
\end{equation}
 for all $x_1,\ldots , x_n\in B(R')$ such that $|x_i|_{\real^d} \geq R$.
\end{enumerate}
 We also let ${\cal S}$ denote the union of the 
 sets ${\cal S}_A$ over the closed subsets $A$ of $B(R')$.
\end{definition}
The gradient operator $D$ is defined on random functionals
 $F \in {\cal S}$ of the form \eqref{f} as 
\begin{equation} 
\label{ndfgd} 
D_y F := 
  \sum_{n=1}^\infty 
   {\bf 1}_{\{ \gamma ( B(R))=n\} }
 \sum_{i=1}^n 
 \langle \G_\eta (X_i,y) ,
 \nabla^{\real^d}_{x_i} f\left( X_1,\ldots , X_n \right)
 \rangle_{\real^d} 
, 
\end{equation} 
$y\in B(R)$. 
For any $F\in {\cal S}$, by \eqref{gxy1} 
we have $DF\in L^1(\Omega \times B(R))$ from the bound 
\begin{eqnarray*} 
  E \left[ \int_{B(R)} | D_x F | \ \! \lambda (dx) \right]
  & \leq & 
\Vert | \nabla^{\real^d} f|_{\real^d} \Vert_\infty
 E \left[
   \int_{B(R)} \int_{B(R)} | \G_\eta (x,y) |_{\real^d}
  \gamma (dx) \lambda (dy)  
 \right]
 \\
 & = &
\Vert | \nabla^{\real^d} f|_{\real^d} \Vert_\infty
 \int_{B(R)} \int_{B(R)} | \G_\eta (x,y) |_{\real^d} 
   \lambda (dx) 
   \lambda (dy) 
\\
 & = &
    K_d  \Vert | \nabla^{\real^d} f|_{\real^d} \Vert_\infty
 \int_{B(R)} \int_{B(R)} \frac{1}{|x-y|_{\real^d}^{d-1}}
    \lambda ( dx ) \lambda (dy)
  \\
 & \leq &
 K_d v_d^2 R' R^d  
\Vert | \nabla^{\real^d} f|_{\real^d} \Vert_\infty
 \\
 & < & \infty. 
\end{eqnarray*}  
\subsubsection*{Poisson-Skorohod integral}
 We let ${\cal U}_0$ denote the space of simple random fields of the form 
\begin{equation}
\label{u} 
 u = \sum_{i=1}^{n} g_i G_i, \qquad n\geq 1, 
\end{equation}
with $G_i \in {\cal S}_{A_i}$ 
and $g_i \in {\cal C}^\infty_0 (B(R))$, 
 $i=1,\ldots ,n$.
  \begin{definition} 
 We define the Poisson-Skorohod integral $\delta (u)$ of $u \in {\cal U}_0$ 
 of the form \eqref{u} as 
\begin{equation}
\label{dltu} 
 \delta(u) := \sum_{i=1}^n
 \Big( G_i \int_{B(R)} g_i (x) ( \gamma (dx) - \lambda (dx) )
 - \langle g_i , DG_i \rangle_{L^2(B(R))} \Big). 
\end{equation}
\end{definition} 
 In particular, for $h\in {\cal C}^\infty_0 (B(R))$
  we have 
$$
 \delta (h) = \int_{B(R)} h(x) ( \gamma (dx) - \lambda (dx) ). 
$$ 
 The proof of the next proposition, cf. 
 Proposition~8.5.1 in \cite{pratprivault}
 and Proposition~5.1 in \cite{prconnection},
 is given in the appendix.
\begin{prop} 
  \label{prop1}
  The operators $D$ and $\delta$ satisfy the duality relation   
\begin{equation} 
  \label{jldf}
  E[ \langle u , DF\rangle_{L^2(B(R))} ] 
 =
 E [ F \delta ( u ) ],
\qquad F\in {\cal S}, 
 \qquad u\in {\cal U}_0. 
\end{equation}
\end{prop} 
As a consequence of Proposition~\ref{prop1}
and the denseness of ${\cal S}$ in $L^1(\Omega )$
and that of ${\cal U}_0$ in $L^1 (\Omega \times B(R) )$,
the gradient operator $D$ is 
 closable in the sense that if $(F_n)_{n\in\inte}\subset {\cal S}$ 
 tends to zero in $L^2 (\Omega )$ 
 and $(DF_n)_{n\in\inte}$ converges 
 to $U$ in $L^1(\Omega \times B(R))$, then $U=0$ a.e.. 
 Similarly, the divergence operator $\delta$ is 
 closable in the sense that if $(u_n)_{n\in\inte}\subset {\cal U}_0$ 
 tends to zero in $L^2 (\Omega \times B(R))$ 
 and $(\delta ( u_n ))_{n\in\inte}$ converges 
 to $G$ in $L^1 (\Omega )$, then $G=0$ a.s.. 
\\

 The gradient operator $D$ defines the Sobolev space 
 $\dee^{1,1}$ with the Sobolev norm 
$$ 
 \Vert F \Vert_{\dee_{1,1} }
 := 
 \Vert F \Vert_{L^2 ( \Omega )}
 + 
 \Vert D F \Vert_{L^1 ( \Omega  \times B(R) )}
, 
 \qquad 
 F \in {\cal S}. 
$$ 
 In the sequel we fix a total order $\preceq$ on $B(R)$ and
 consider the space ${\cal P}_0\subset {\cal U}_0$ 
  of simple predictable random field of the form 
\begin{equation}
\label{uf} 
u: = \sum_{i=1}^n g_iF_i, 
\end{equation} 
such that the supports of $g_1,\ldots , g_n$ satisfy 
$$
\Supp ( g_i ) \preceq \cdots \preceq \Supp ( g_n )
\quad
\mbox{and}
\quad
F_i \in {\cal S}_{A_i}, 
$$ 
where $\Supp ( g_1 )\cup \cdots \cup \Supp ( g_{i-1} ) \subset A_i\subset B(R')$
and $A_i \preceq \Supp ( g_i )$, $i=1,\ldots , n$.
\\ 

Such random fields are predictable in the sense of e.g. \S~5 of
\cite{last} and references therein.
\\

We will also assume that the order $\preceq$ is compatible with the
kernel $\G_\eta $ in the sense that 
\begin{equation}
  \label{comp}
  \G_\eta (x,y) = 0
  \quad
  \mbox{for all} \quad
  x, y \in B(R)
  \quad
  \mbox{such that} \quad
  x \preceq y.
\end{equation}
 Under the compatibility condition \eqref{comp} we have in particular
$$
D_y F = 0, \quad
y\in B(R),
\quad
A\preceq y,
\qquad
F\in {\cal S}_A.
$$ 
Moreover, if $u \in {\cal P}_0$ is a predictable random field
of the form \eqref{uf} we note that by \eqref{ndfgd} and the compatibility condition
\eqref{comp} we have 
$$
D_yF_i=0, \qquad A_i \preceq y, \quad i=1,\ldots , n,
$$
 hence 
\begin{equation}
\label{comp2} 
D_y u_x = 0 ,
\qquad
x \preceq y,
\quad
x,y \in B(R).
\end{equation}
{\em Example}.
The order $\preceq$ defined by
\begin{equation}
  \label{compat} 
x = (x^{(1)},\ldots , x^{(d)}) \preceq y = (y^{(1)},\ldots , y^{(d)})
\quad
\Longleftrightarrow
\quad
x^{(1)} \leq y^{(1)}
\end{equation} 
is compatible with the
kernel $\G_\eta$ provided that the support of
$\eta$ is contained in
$$
\big\{
x = (x^{(1)},\ldots , x^{(d)}) \in B(R')\setminus B(R) \
: \
x^{(1)}>R \big\}.
$$
 
The proof of the next Proposition~\ref{prop2}
is given in the appendix.
\begin{prop} 
\label{prop2}
 The Poisson-Skorohod integral of 
 $u=(u_x)_{x\in B( R )}$ in the space ${\cal P}_0$ of simple predictable
 random fields satisfies the relation 
\begin{equation}
  \label{rel} 
\delta ( u ) = \int_{B(R)} u_x ( \gamma ( dx ) - \lambda ( dx ) ), 
\end{equation} 
which extends to the closure of ${\cal P}_0$ in $L^2(\Omega \times B(R))$
by density and the isometry relation 
\begin{equation}
  \label{isomrel} 
E [ \delta (u)^2 ]
=
E\left[
\int_{B(R)} u^2_x \ \! \lambda ( dx ) 
  \right], \qquad u \in {\cal P}_0.
\end{equation} 
\end{prop} 
\subsubsection*{Covariant derivative} 
 In addition to the gradient operator $D$,
 we will also need the following notion of covariant
 derivative operator $\widetilde{\nabla}$ defined on
 stochastic processes that are viewed as tangent processes on the Poisson space
 $\Omega$, see \cite{prconnection}. 
\begin{definition} 
  Let the operator $\widetilde{\nabla}$ be defined
  on $u \in {\cal P}_0$ as 
\begin{equation} 
\nonumber 
 \widetilde{\nabla}_y u_x : =
 {D}_y u_x + \langle \G_\eta (x,y) , \nabla^{\real^d}_x u_x \rangle_{\real^d}, 
 \qquad x, y \in B(R). 
\end{equation} 
\end{definition} 
We note that from the compatibility condition \eqref{comp}
and Relation~\eqref{comp2} we also have
\begin{equation}
\label{comp3} 
\widetilde{\nabla}_y u_x = 0 ,
\qquad
x \preceq y,
\quad
x,y \in B(R).
\end{equation}
From the bound
\begin{eqnarray*} 
  \lefteqn{
    \!  \!  \!  \!  \!  \!  \!  \!  \!  \!  \!
    E\left[ \int_{B(R)\times B(R)}| \widetilde{\nabla}_x u_y | \ \! \lambda ( dx ) \lambda (dy) \right] 
  }
  \\
  & \leq &  
\Vert D u \Vert_{L^1(\Omega \times B(R) \times B(R) )}
  +
  E\left[
    \int_{B(R) \times B(R) } | \langle \G_\eta (x,y) , \nabla^{\real^d}_x u_x \rangle_{\real^d} |
    \ \! \lambda ( dx ) \lambda (dy)
    \right]
  \\
  & \leq &
 \Vert D u \Vert_{L^1(\Omega \times B(R) \times B(R) )}
  +
  K_d E\left[
    \int_{B(R) \times B(R) } \frac{1}{|x-y|_{\real^d}^{d-1}}
    | \nabla^{\real^d} u_x |_{\real^d} 
    \lambda ( dx ) \lambda (dy)
    \right]
  \\
  & \leq &
 \Vert D u \Vert_{L^1(\Omega \times B(R) \times B(R) )}
  + K_d v_d R'
  E\left[
    \int_{B(R)}
    | \nabla^{\real^d}_x u_x |_{\real^d} 
    \lambda ( dx ) 
    \right]
  \\
  & = &
 \Vert D u \Vert_{L^1(\Omega \times B(R) \times B(R) )}
  + K_d v_d R' \Vert \nabla^{\real^d} u \Vert_{L^1(\Omega \times B(R) ; \real^d)}, 
\end{eqnarray*} 
we check that $\widetilde{\nabla}$ extends to the Sobolev space
$\widetilde{\dee}^{1,1}_0$ of {\em predictable} random fields defined
as the completion of
${\cal P}_0$ under the Sobolev norm  
$$ 
 \Vert u \Vert_{ \widetilde{\dee}^{1,1}} 
 : = 
 \Vert u \Vert_{L^2(\Omega , W^{1,1}_0(B(R)))} 
 + \Vert D u \Vert_{L^1(\Omega \times B(R)\times B(R) ) } 
 ,
 \qquad u \in {\cal P}_0,
$$ 
 where $W^{1,p}_0 (B(R))$ is the first order Sobolev space completion
 of ${\cal C}^\infty_0(B(R))$ under the norm 
$$ 
 \Vert f \Vert_{W^{1,p}(B(R))} 
 : = 
 \Vert f \Vert_{L^p(B(R))} 
 + 
 \Vert \nabla^{\real^d} f \Vert_{L^p(B(R);\real^d)}, \qquad p \geq 1. 
$$
\subsubsection*{Commutation relation} 
 In the sequel, we denote by $\widetilde{\dee}^{1,\infty}_0$ the set of
 predictable random fields $u$ in $\widetilde{\dee}^{1,1}_0$
 that are bounded together with their covariant derivative 
 $\widetilde{\nabla}u$. 
\begin{prop} 
\label{commut} 
For $u \in \widetilde{\dee}^{1,\infty}_0$ a predictable random field,
we have the commutation relation 
\begin{equation} 
\label{com} 
D_y \delta(u) = u(y) + \delta(\widetilde{\nabla}_y u),
\qquad y\in B(R). 
\end{equation}
\end{prop} 
\begin{Proof} 
 Taking $h\in {\cal C}^\infty_0 (B(R))$, we have $\delta( h) \in {\cal S}$ and
\begin{eqnarray*} 
 D_y \delta( h) & = & 
 D_y \int_{B(R)} h(y) ( \gamma (dx ) - \lambda ( dx ) ) 
\\
 & = & 
 \int_{B(R)} \langle \G_\eta (x,y) , \nabla^{\real^d}_x h(x) \rangle_{\real^d} \gamma (dx) 
\\
 & = &
  \int_{B(R)} \langle \G_\eta (x,y) , \nabla^{\real^d}_x h(x) \rangle_{\real^d} \lambda (dx)
 +
 \delta ( \widetilde{\nabla}_y h ) 
 \\
 & = &
 h(y) + \delta ( \widetilde{\nabla}_y h ). 
\end{eqnarray*} 
 where we applied \eqref{gxy4-2}. 
 Next, taking $u=hF \in {\cal P}_0$ a simple predictable random field,
 we check that $\delta( u) \in {\cal S}$, and by \eqref{dltu} or \eqref{dltu2} we have 
\begin{eqnarray*} 
 D_y \delta( F h) & = & 
 D_y \left( F \delta( h) - \langle h , DF\rangle_{L^2(B(R))} \right)
 \\
 & = &
 D_y \left( F \delta( h) \right)
 \\
 & = &
 \delta( h) D_y F + F D_y \delta( h) 
 \\
 & = &
 \delta( h) D_y F + F ( h(y) + \delta ( \widetilde{\nabla}_y h ) )
 \\
 & = &
 F h(y) + \delta( h D_y F + F \widetilde{\nabla}_y h )
 \\
 & = &
 F h(y)
 +
  \delta( \widetilde{\nabla}_y ( F h ) ) 
\\
 & = &
 u_y + \delta( \widetilde{\nabla}_y u ) ,
 \qquad y \in B(R). 
\end{eqnarray*} 
We conclude by the denseness of ${\cal P}_0$ in $\widetilde{\dee}^{1,1}_0$
and by the closability of $\widetilde{\nabla}$, $D$ and $\delta$.
\end{Proof}
\section{Cumulant operators}
\label{s2-1}
In the sequel, given $h$ in the standard Sobolev space
$W^{1,p}( B(R) )$ on $B(R)$ and $f \in L^q(B(R))$ 
with $1=p^{-1}+q^{-1}$, $p,q\in [1,\infty]$, we define 
\begin{equation} 
\label{1lpa} 
 ( \widetilde{\nabla} h ) f_x
: =
 \int_{B(R)} f(y) \widetilde{\nabla}_y h(x) \ \! \lambda (dy)
 =
 \int_{B(R)} f(y) \langle \G_\eta (x,y) , \nabla^{\real^d}_x h(x) \rangle_{\real^d} \lambda (dy)
, 
\end{equation} 
$x \in B(R)$. 
More generally, given $k\geq 1$ and $u \in \widetilde{\dee}^{1,1}_0$ 
a predictable random field, we let the operator $( \widetilde{\nabla} u)^k$ 
be defined in the sense of matrix powers with continuous indices,
as 
\begin{equation} 
\nonumber 
 ( \widetilde{\nabla} u )^k 
 f_y 
 = 
 \int_{B(R)}
 \cdots 
 \int_{B(R)}
 ( 
 \widetilde{\nabla}_{x_k} u_y 
 \widetilde{\nabla}_{x_{k-1}} u_{x_k} 
 \cdots 
 \widetilde{\nabla}_{x_1} u_{x_2} 
 ) 
 f_{x_1} 
 \ \! \lambda ( dx_1 ) \cdots \lambda ( dx_k), 
\end{equation} 
$y \in B(R)$, $f \in L^2(B(R))$.
\begin{prop}
  \label{p1} 
  For any $n \in \inte$, $p > 1$, $r\in [0,1]$, 
  $h\in W^{1,p/(1-r)^{n-1}/r}( B(R) )$ and $f \in L^{p/(1-r)^n} (B(R))$ we have the bound
\begin{equation}
\label{vrt} 
  \Vert ( \widetilde{\nabla} h )^n f \Vert_{L^p(B(R))}
 \leq 
  ( K_d v_d R' )^n 
  \Vert f \Vert_{L^{p/(1-r)^n}(B(R))}
  \prod_{j=1}^n \Vert \nabla^{\real^d} h \Vert_{L^{p/(1-r)^{j-1}/r}(B(R);\real^d)}. 
\end{equation} 
\end{prop} 
\begin{Proof}
 For $n=1$ we have 
\begin{eqnarray} 
  \nonumber
  \lefteqn{
    \Vert ( \widetilde{\nabla} h ) f \Vert_{L^p(B(R))}^p
  = 
  \int_{B(R)} \left| \int_{B(R)} f(y) \widetilde{\nabla}_y h(x) \ \! \lambda (dy) \right|^p
  \lambda (dx)
  }
  \\
\nonumber
    & = &
  \int_{B(R)} \left| \int_{B(R)} f(y)
  \langle \G_\eta (x,y) , \nabla^{\real^d}_x h(x) \rangle_{\real^d}
  \lambda (dy) \right|^p
  \lambda (dx)
  \\
\nonumber
    & = &
  \int_{B(R)} \left| \left< 
  \int_{B(R)} f(y)
  \G_\eta (x,y) \ \! \lambda (dy)
  , \nabla^{\real^d}_x h(x) \right>_{\real^d}
  \right|^p
  \lambda (dx)
  \\
\nonumber
    & \leq &
  \int_{B(R)} 
  \left| 
  \int_{B(R)} f(y)
  \G_\eta (x,y) \ \! \lambda (dy)
  \right|_{\real^d}^p
  | \nabla^{\real^d}_x h(x) |_{\real^d}^p 
  \lambda (dx)
  \\
\nonumber
    & = &
  \left( 
    \int_{B(R)} 
  \left| 
  \int_{B(R)} f(y)
  \G_\eta (x,y) \ \! \lambda (dy)
  \right|_{\real^d}^{p/(1-r)}
  \lambda (dx)
  \right)^{1-r} 
  \left( 
    \int_{B(R)} 
  | \nabla^{\real^d}_x h(x) |_{\real^d}^{p/r} 
  \lambda (dx)
 \right)^r 
  \\
  \label{djcdfd}
  & \leq &
  ( K_d v_d R' )^p 
  \Vert f \Vert_{L^{p/(1-r)}(B(R))}^p
  \Vert \nabla^{\real^d} h \Vert_{L^{p/r}(B(R);\real^d)}^p,
  \end{eqnarray} 
 where we used the bound \eqref{gxy3}. 
 Next, assuming that \eqref{vrt} holds at the rank $n\geq 1$
 and using \eqref{djcdfd}, we have
\begin{eqnarray*} 
  \Vert ( \widetilde{\nabla} h )^{n+1} f \Vert_{L^p(B(R))}
  & = &
  \Vert ( \widetilde{\nabla} h )^n ( \widetilde{\nabla} h ) f \Vert_{L^p(B(R))}
  \\
   & \leq &  
  ( K_d v_d R' )^n 
  \Vert ( \widetilde{\nabla} h ) f \Vert_{L^{p/(1-r)^n}(B(R))}
  \prod_{j=1}^n \Vert \nabla^{\real^d} h \Vert_{L^{p/(1-r)^{j-1}/r}(B(R);\real^d)}
    \\
   & \leq &  
  ( K_d v_d R' )^{n+1}
  \Vert f \Vert_{L^{p/(1-r)^{n+1}} (B(R))}
  \prod_{j=1}^{n+1} \Vert \nabla^{\real^d} h \Vert_{L^{p/(1-r)^{j-1}/r}(B(R);\real^d)}, 
\end{eqnarray*} 
 and we conclude to \eqref{vrt} by induction. 
\end{Proof}
In particular, for $r=0$, $f\in L^p (B(R))$, $p>1$,
and $h\in W^{1,1}(B(R))$ the argument of Proposition~\ref{p1}
shows that 
$$ 
 \Vert ( \widetilde{\nabla} h )^n f \Vert_{L^p(B(R))}
 \leq 
  ( K_d v_d R' )^n 
  \Vert f \Vert_{L^p (B(R))}
  \Vert \nabla^{\real^d} h \Vert_{L^\infty (B(R);\real^d)}^n,
  \qquad  n \in \inte. 
$$ 
 We note that for $u \in \widetilde{\dee}^{1,\infty}_0$ a  
 predictable random field, the
 random field $( \widetilde{\nabla} u) u \in \widetilde{\dee}^{1,\infty}_0$
 is also predictable from \eqref{comp3} and \eqref{1lpa}.
\\

 In the next definition we construct a family of cumulant operators
    which differs from the one introduced in \cite{nourdin10} on the Wiener space. 
\begin{definition}
\label{def01}
Given $k\geq 2$ and $u \in \widetilde{\dee}^{1,\infty}_0$ a  
predictable random field 
we define the operators 
$\Gamma^u_k : \dee_{1,1} \longrightarrow L^1( \Omega )$ by 
$$ 
 \Gamma^u_k F 
 : = 
 F 
 \langle 
 ( \widetilde{\nabla} u )^{k-2} u 
 , 
 u 
 \rangle_{L^2(B(R))} 
 + 
 \langle 
 ( \widetilde{\nabla} u )^{k-1} u 
 , 
 DF 
 \rangle_{L^2(B(R))}
 , 
 \quad
F\in \dee_{1,1}. 
$$ 
\end{definition}
  We note that for 
  $h$ in the space $W^{1,\infty}( B(R) )$ of bounded functions in 
  $W^{1,1}( B(R) )$, and $f \in L^p (B(R))$, $p>1$,
  $m\geq 1$, we have 
\begin{eqnarray*} 
 \langle h^m , ( \widetilde{\nabla} h ) f \rangle_{L^2(B(R))} 
 & = &
 \int_{B(R)} h^m(x)  \int_{B(R)} f(y) \langle \G_\eta (x,y)  , \nabla^{\real^d}_x h(x) \rangle_{\real^d}
 \lambda (dy)
 \lambda (dx)
 \\
  & = &
 \frac{1}{m+1}
 \int_{B(R)} \int_{B(R)} f(y) \langle \G_\eta (x,y)  , \nabla^{\real^d}_x h^{m+1} (x) \rangle_{\real^d}
 \lambda (dy)
 \lambda (dx)
 \\
  & = &
 \frac{1}{m+1}
 \int_{B(R)}
 f(x) h^{m+1} (x) 
 \ \! \lambda (dx), 
\end{eqnarray*}
where we applied \eqref{gxy4}, hence 
$$ 
 \langle h^m , ( \widetilde{\nabla} h )^{n+1} f \rangle_{L^2(B(R))} 
 =
 \frac{1}{m+1}
 \int_{B(R)} 
 h^{m+1} (x) ( \widetilde{\nabla} h )^n f (x)  
 \ \! \lambda (dx)
 ,
 $$
 which implies by induction
 $$ 
 \langle ( \widetilde{\nabla} h )^n f , h^m \rangle_{L^2(B(R))} 
 =
 \frac{m!}{(m+n)!}
 \int_{B(R)} 
 h^{m+n} (x) f (x)  
 \ \! \lambda (dx). 
$$ 
 In Lemma~\ref{itlw2.1} we generalize this identity to $h$ a random field. 
\begin{lemma} 
\label{itlw2.1} 
For $n\in \inte$, $m \geq 1$, $u \in \widetilde{\dee}^{1,\infty}_0$
a predictable random field and $f \in L^p (B(R))$, $p>1$, we have 
\begin{eqnarray} 
\label{djklll} 
 \langle 
 ( \widetilde{\nabla} u )^n f  
 , 
 u^m 
 \rangle_{L^2(B(R))} 
 & = & 
 \frac{m!}{(m+n)!} 
 \int_{B(R)}
 u_x^{m+n} f (x) \ \! \lambda ( dx )  
\\
 \nonumber
 & &
 + \sum_{k=1}^n 
 \frac{m!}{(m+k)!} 
 \left< 
 ( \widetilde{\nabla} u )^{n-k} 
 f  , 
 D 
 \int_{B(R)}
 u^{m+k}_x \ \! \lambda ( dx )  
 \right>_{L^2(B(R))}.
\end{eqnarray} 
\end{lemma} 
\begin{Proof} 
  Using the adjoint $\widetilde{\nabla}^* u$ of $\widetilde{\nabla} u$ on $L^2(B(R))$
  given by 
$$ 
 ( \widetilde{\nabla}^* u) v_y : = \int_{B(R)} ( \widetilde{\nabla}_y  u_x ) v_x \ \! \lambda ( dx ), 
 \qquad 
 y \in B(R), \quad 
 v \in L^2(B(R)), 
$$
 with the duality relation 
$$ 
 \langle v , ( \widetilde{\nabla}^* u ) h \rangle_{L^2(B(R))}
 =
  \langle ( \widetilde{\nabla} u ) v , h \rangle_{L^2(B(R))}, 
 \qquad 
 h ,v \in L^2(B(R)), 
$$ 
 we will show by induction on $k=0,1,\ldots, n$ that 
\begin{align} 
\nonumber 
& ( \widetilde{\nabla}^* u )^n u^m_{x_0} 
 = 
 \int_{B(R)} 
 \cdots 
 \int_{B(R)} 
 u^m_{x_n} 
 \widetilde{\nabla}_{x_0} u_{x_1} 
 \widetilde{\nabla}_{x_1} u_{x_2} 
 \cdots 
 \widetilde{\nabla}_{x_{n-1}} u_{x_n} 
 \ \! \lambda ( dx_1 ) \cdots \lambda ( dx_n )  
\\
 \nonumber
 & = 
\sum_{i=1}^k 
 \frac{m!}{(m+i)!}
 \int_{B(R)} 
 \cdots 
 \int_{B(R)} 
 \widetilde{\nabla}_{x_0} u_{x_1} 
 \cdots 
 \widetilde{\nabla}_{x_{n-i-1}} u_{x_{n-i}} 
 D_{x_{n-i}} 
 u_{x_{n+1-i}}^{m+i}
 \ \! \lambda ( dx_1 ) \cdots \lambda ( dx_{n-i-1} ) 
\\ 
\label{id} 
 & 
 + 
 \frac{m!}{(m+k)!} 
 \int_{B(R)} 
 \cdots 
 \int_{B(R)} 
 u_{x_{n-k}}^{m+k} 
 \widetilde{\nabla}_{x_0} u_{x_1} 
 \cdots 
 \widetilde{\nabla}_{x_{n-k-1}} u_{x_{n-k}} 
 \ \! \lambda ( dx_1 ) \cdots \lambda ( dx_{n-k} ). 
\end{align} 
 By \eqref{1lpa}, this relation holds for $k=0$. 
 Next, assuming that the identity \eqref{id} holds for some 
 $k \in \{ 0,1, \ldots ,n-1 \}$, and using the relation 
$$ 
 \widetilde{\nabla}_{x_{n-k-1}} u_{x_{n-k}} 
 = 
 D_{x_{n-k-1}} u_{x_{n-k}} 
 + 
 \langle  \G_\eta (x_{n-k}, x_{n-k-1} ) , \widetilde{\nabla}_{x_{n-k}} u_{x_{n-k}} \rangle_{\real^d},
 \quad 
 x_{n-k-1} , x_{n-k} \in B(R), 
$$ 
 we have 
\begin{align*} 
& ( \widetilde{\nabla}^* u )^n u_{x_0} 
  \\
  &
  =  
 \sum_{i=1}^k 
 \frac{m!}{(m+i)!} 
 \int_{B(R)} 
 \cdots 
 \int_{B(R)}
 \widetilde{\nabla}_{x_0} u_{x_1} 
 \cdots 
 \widetilde{\nabla}_{x_{n-i-1}} u_{x_{n-i}} 
 D_{x_{n-i}} 
 u_{x_{n+1-i}}^{m+i} 
 \ \! \lambda ( dx_1 ) \cdots \lambda ( x_{n+1-i} ) 
\\ 
 &  
 + 
 \frac{m!}{(m+k)!} 
 \int_{B(R)} 
 \cdots 
 \int_{B(R)} 
 u_{x_{n-k}}^{m+k} 
 \widetilde{\nabla}_{x_0} u_{x_1} 
 \cdots 
 \widetilde{\nabla}_{x_{n-k-1}} u_{x_{n-k}} 
 \ \! \lambda ( dx_1 ) \cdots \lambda ( dx_{n-k} )  
\\ 
 & =  
 \sum_{i=1}^k 
 \frac{m!}{(m+i)!} 
 \int_{B(R)} 
 \cdots 
 \int_{B(R)} 
 \widetilde{\nabla}_{x_0} u_{x_1} 
 \cdots 
 \widetilde{\nabla}_{x_{n-i-1}} u_{x_{n-i}} 
 D_{x_{n-i}} 
 u_{x_{n+1-i}}^{m+i} 
 \ \! \lambda ( dx_1 ) \cdots \lambda ( dx_{n+1-i} ) 
\\ 
 &  
 + 
 \frac{m!}{(m+k)!} 
 \int_{B(R)} 
 \cdots 
 \int_{B(R)} 
 u_{x_{n-k}}^{m+k} 
 \widetilde{\nabla}_{x_0} u_{x_1} 
 \cdots 
 \widetilde{\nabla}_{x_{n-k-2}} u_{x_{n-k-1}} 
 D_{x_{n-k-1}} u_{x_{n-k}} 
 \ \! \lambda ( dx_1 ) \cdots \lambda ( dx_{n-k} ) 
\\ 
 &  
 + 
 \frac{m!}{(m+k)!} 
 \int_{B(R)} 
 \cdots 
 \int_{B(R)}
 \langle  \G_\eta (x_{n-k} , x_{n-k-1} ) , \widetilde{\nabla}_{x_{n-k}} u_{x_{n-k}} \rangle_{\real^d}
 \\
 & \qquad
 \times 
 u_{x_{n-k}}^{m+k-2} 
 \widetilde{\nabla}_{x_0} u_{x_1} 
 \cdots 
 \widetilde{\nabla}_{x_{n-2-k}} u_{x_{n-k-1}} 
 \ \! \lambda ( dx_1 ) \cdots \lambda ( dx_{n-k} ) 
\\ 
 & =  
 \sum_{i=1}^k 
 \frac{m!}{(m+i)!} 
 \int_{B(R)} 
 \cdots 
 \int_{B(R)} 
 \widetilde{\nabla}_{x_0} u_{x_1} 
 \cdots 
 \widetilde{\nabla}_{x_{n-i-1}} u_{x_{n-i}} 
 D_{x_{n-i}} 
 u_{x_{n+1-i}}^{m+i} 
 \ \! \lambda ( dx_1 ) \cdots \lambda ( dx_{n+1-i} ) 
\\ 
 &  
 + 
 \frac{m!}{(m+k+1)!} 
 \int_{B(R)} 
 \cdots 
 \int_{B(R)}
 \widetilde{\nabla}_{x_0} u_{x_1} 
 \cdots 
 \widetilde{\nabla}_{x_{n-k}} u_{x_{n-k-1}} 
 D_{x_{n-k-1}} u_{x_{n-k}}^{m+k+1} 
 \ \! \lambda ( dx_1 ) \cdots \lambda ( dx_{n-k} ) 
\\ 
 &  
 + 
 \frac{m!}{(m+k+1)!} 
 \int_{B(R)} 
 \cdots 
 \int_{B(R)} 
 \widetilde{\nabla}_{x_0} u_{x_1} 
 \cdots 
 \widetilde{\nabla}_{x_{n-k-2}} u_{x_{n-k-1}} 
 \\
 & \qquad
 \times \int_{B(R)} 
 \langle  \G_\eta (x , x_{n-k-1} ) , \nabla^{\real^d}_x u^{m+k+1}_x \rangle_{\real^d}
 \lambda ( dx )  
 \lambda ( dx_1 ) \cdots \lambda ( dx_{n-k-1} ) 
\\ 
 & =  
 \sum_{i=1}^{k+1} 
 \frac{m!}{(m+i)!} 
 \int_{B(R)} 
 \cdots 
 \int_{B(R)} 
 \widetilde{\nabla}_{x_0} u_{x_1} 
 \cdots 
 \widetilde{\nabla}_{x_{n-i-1}} u_{x_{n-i}} 
 D_{x_{n-i}} 
 u_{x_{n+1-i}}^{m+i} 
 \ \! \lambda ( dx_1 ) \cdots \lambda ( dx_{n+1-i} ) 
\\ 
 &  
 + 
 \frac{m!}{(m+k+1)!} 
 \int_{B(R)} 
 \cdots 
 \int_{B(R)} 
 u_{x_{n-k-1}}^{m+k+1} 
 \widetilde{\nabla}_{x_0} u_{x_1} 
 \cdots 
 \widetilde{\nabla}_{x_{n-k-2}} u_{x_{n-k-1}} 
 \ \! \lambda ( dx_1 ) \cdots \lambda ( dx_{n-k-1} ) 
\\ 
 & =  
 \sum_{i=1}^{k+1} 
 \frac{m!}{(m+i)!} 
 ( \widetilde{\nabla}^* u )^{n-i} 
 D_{x_0} 
 \int_{B(R)}
 u_s^{m+i} 
 \ \! \lambda ( ds )  
 + 
 \frac{m!}{(m+k+1)!} 
 ( \widetilde{\nabla}^* u )^{n-k-1} 
 u_{x_0}^{m+k+1} 
, 
\end{align*} 
which shows by induction that \eqref{id} holds
at the rank $k=n$, in particular we have 
\begin{equation} 
\nonumber 
 ( \widetilde{\nabla}^* u )^n u^m_x 
 = 
 \frac{m!}{(m+k)!} 
 u^{m+n}_x 
 + 
 \sum_{i=2}^{n+1} 
 \frac{m!}{(m+i-1)!} 
 ( \widetilde{\nabla}^* u )^{n+1-i} 
 D_x 
 \int_{B(R)}
 u_y^{m+i-1} 
 \ \! \lambda ( dy ) 
, 
\end{equation} 
$x \in B(R)$, which yields \eqref{djklll} 
by integration with respect to $x \in B(R)$ and duality. 
\end{Proof} 
 As a consequence of Lemma~\ref{itlw2.1} we have
\begin{equation} 
\nonumber 
 \Gamma^u_k {\bf 1}  
 = 
 \int_{B(R)}
 \frac{u_x^k }{(k-1)!} \ \! \lambda ( dx ) 
 + 
 \sum_{i=2}^{k-1} 
 \frac{1}{i!} 
 \left< 
 ( \widetilde{\nabla} u )^{k-1-i} 
 u , 
 D 
 \int_{B(R)} 
 u^i_x \ \! \lambda ( dx ) 
 \right>_{L^2(B(R))},   
\end{equation} 
$k\geq 2$.  
Hence when $h\in W^{1,p}( B(R) )$, $p>1$,
is a deterministic 
function such that $\Vert \nabla^{\real^d} h\Vert_\infty < \infty$, 
 we find the relation 
\begin{equation} 
\label{ifinad2} 
 \Gamma^h_k {\bf 1} 
 = 
 \frac{1}{(k-1)!} 
 \int_{B(R)} 
 h^k(x) \ \! \lambda ( dx ) 
 = 
 \frac{1}{(k-1)!} 
 \kappa_k^h 
 , \qquad
 k \geq 2,
\end{equation} 
 which shows that $\Gamma^h_k {\bf 1}$ coincides with 
 the cumulant $\kappa_k^h = \int_{B(R)} h^k(x) \ \! \lambda ( dx )$
 of order $k\geq 2$ of the Poisson stochastic
 integral $\int_{B(R)} h(x) ( \gamma ( dx ) - \lambda (dx) )$. 
\section{Edgeworth type expansions} 
\label{s3} 
 Classical Edgeworth series provide expansion 
 of the cumulative distribution function $P(F\leq x)$ 
 of a centered random variable $F$ with $E[F^2]=1$
 around the Gaussian cumulative distribution function 
 $\Phi (x)$, using the cumulants $(\kappa_n)_{n\geq 1}$ of a random variable $F$
 and Hermite polynomials. 
 Edgeworth type expansions of the form 
\begin{equation} 
\nonumber 
 E[Fg(F) ] = \sum_{l=1}^n \frac{\kappa_{l+1}}{l!} E [ g^{(l)} (F) ] 
 + E[g^{(n+1)} (F) \Gamma_{n+1} F ], \quad n \geq 1, 
\end{equation} 
 for $F$ a centered random variable, 
 have been obtained by the Malliavin calculus in 
 \cite{nourdinpeccati2009}, 
 where $\Gamma_{n+1}$ is a cumulant type 
 operator on the Wiener space such that 
 $n! E[\Gamma_n F ]$ coincides with the cumulant
 $\kappa_{n+1}$ of order $n+1$ of $F$, $n\in \inte$, cf. 
 \cite{nourdin10}, extending the results of \cite{barbour1986}
 to the Wiener space.
\\
 
In this section we establish an Edgeworth type expansion
of any finite order with an explicit remainder term for the
compensated Poisson stochastic integral $\delta (u)$ of
a 
predictable random field $(u_x)_{x\in B(R)}$.
In the sequel we let $\langle \cdot , \cdot \rangle$
denote $\langle \cdot , \cdot \rangle_{L^2(B(R))}$, except
if stated otherwise. 
\\
 
 Before proceeding to the statement of
 general expansions in Proposition~\ref{djklsad1}, 
 we illustrate the method with the derivation of an expansion of order one for
 a deterministic integrand $f$. 
 By the duality relation \eqref{jldf} between $D$ and $\delta$, 
 the chain rule of derivation for $D$ and 
 the commutation relation \eqref{com} we get, 
 for ${{g}} \in {\cal C}^2_b(\real )$ and 
 $f\in W^{1,1}_0(B(R))$ such that
 $\Vert \nabla^{\real^d} f\Vert_\infty < \infty$, 
\begin{eqnarray*} 
  \lefteqn{
    E [ 
 \delta ( f ) 
 {{g}} ( \delta ( f ) ) 
 ] 
 = 
 E [ 
 \langle f , D \delta ( f ) \rangle 
         {{g}}' ( \delta ( f ) )
   ] 
  }
  \\ 
 & = & 
 E[ 
 \langle 
 f , 
 f 
 \rangle 
 g' ( \delta ( f ) )
 ] 
 + 
 E \big[
   \langle f , \delta ( \widetilde{\nabla}^* f) \rangle 
   g' ( \delta (f) )
 \big] 
 \\ 
 & = & 
 E[ 
 \langle 
 f , 
 f 
 \rangle 
 g' ( \delta ( f ) )
 ] 
+ 
 E \big[ 
 \langle \widetilde{\nabla}^* f , 
 D ( 
 g' ( \delta (f) ) 
 f 
 ) 
 \rangle 
 \big] 
 \\ 
 & = &  
 E[ 
 \langle 
 f , 
 f 
 \rangle 
 g' ( \delta ( f ) ) 
 ] 
 +
  E\big[ 
 \langle 
 ( \widetilde{\nabla} f ) 
 f  
 , 
 D 
 \delta ( f ) 
 \rangle 
 g'' ( \delta ( f ) ) 
 \big] 
\\ 
 & = & 
 E[ 
 \langle 
 f , 
 f 
 \rangle 
 g' ( \delta ( f ) ) 
 ] 
 + 
\frac{1}{2}
\int_{B(R)} f^3(x) \ \! \lambda ( dx ) 
 E[ 
 g'' ( \delta ( f ) ) 
 ] 
 + 
 E\big[ 
 \langle 
 ( \widetilde{\nabla} f )
 f  
 , 
 \delta ( \widetilde{\nabla}^* f ) 
 \rangle 
 g'' ( \delta ( f ) ) 
 \big]
\\ 
 & = & 
\kappa_2^f 
 E[ 
 g' ( \delta ( f ) ) 
 ] 
 + 
\frac{1}{2}
\kappa_3^f 
 E[ 
 g'' ( \delta ( f ) ) 
 ] 
 + 
 E\big[ 
 g'' ( \delta ( f ) ) 
 \delta (
 ( \widetilde{\nabla} f )^2 f  
 )
 \big]
, 
\end{eqnarray*} 
 since by Lemma~\ref{itlw2.1} we have 
$$ 
\langle (\widetilde{\nabla} f)f , f\rangle 
 = 
\frac{1}{2}
\int_{B(R)} f^3(x) \ \! \lambda ( dx ) 
=
\frac{1}{2}
\kappa_3^f. 
$$ 
In the next proposition we derive general Edgeworth type
expansions for predictable integrand processes $(u_x)_{x \in \real^d}$.
\begin{prop} 
\label{djklsad1} 
Let $u \in \widetilde{\dee}^{1,\infty}_0$ and $n \geq 0$.
For all ${{g}} \in {\cal C}^{n+1}_b(\real )$ 
 and bounded $G\in \dee_{1,1}$ we have 
\begin{align} 
\nonumber 
 &   E\left[ 
 G 
 \delta ( u ) 
 {{g}} ( \delta ( u ) ) 
 \right] 
  = 
 E\left[ 
 \langle u , DG\rangle 
 g ( \delta ( u ) )
 \right] 
 +
 \sum_{k=1}^n 
 E\left[ 
 {{g}}^{(k)} ( \delta ( u ) ) 
 \Gamma_{k+1}^u G 
 \right] 
\\ 
\nonumber 
 & 
+ 
 E\left[ 
 G 
       {{g}}^{(n+1)} ( \delta ( u ) )
 \left(
 \int_{B(R)} 
 \frac{ u_x^{n+2} }{(n+1)!}
 \ \! \lambda ( dx ) 
 + \sum_{k=2}^{n+1} 
 \left< 
 ( \widetilde{\nabla} u )^{n+1-k} 
 u , 
 D \int_{B(R)} \frac{u^k_x}{k!} \ \! \lambda ( dx ) 
 \right> 
 \right)
       \right] 
\\ 
\nonumber 
 & 
+ 
 E\left[ 
 G 
       {{g}}^{(n+1)} ( \delta ( u ) )
 \langle 
 ( \widetilde{\nabla} u )^n 
 u  
 , 
 \delta ( \widetilde{\nabla}^* u ) 
 \rangle 
       \right] 
. 
\end{align} 
\end{prop} 
\begin{Proof} 
 By the duality relation \eqref{jldf} between $D$ and $\delta$, 
 the chain rule of derivation for $D$ and 
 the commutation relation \eqref{com}, we get 
\begin{align*} 
\nonumber 
 & E\big[ 
 G 
 \langle 
 ( \widetilde{\nabla} u )^k 
 u , 
 D 
 \delta ( u )
 \rangle 
 {{g}} ( \delta ( u ) ) 
 \big] 
 - 
 E\big[ 
 G \langle 
 ( \widetilde{\nabla} u )^{k+1} 
 u  
 , 
 D 
 \delta ( u ) 
 \rangle 
 {{g}}' ( \delta ( u ) ) 
 \big] 
\\ 
 & =  
 E\big[ 
 G \langle 
 ( \widetilde{\nabla} u )^k 
 u , 
 u 
 \rangle 
 {{g}} ( \delta ( u ) ) 
 \big] 
 + 
 E \big[ G 
   \langle ( \widetilde{\nabla} u)^k u , \delta ( \widetilde{\nabla}^* u) \rangle 
           {{g}} ( \delta (u) )
 \big] 
  - 
 E\big[ 
 G \langle 
 ( \widetilde{\nabla} u )^{k+1} 
 u  
 , 
 D 
 \delta ( u ) 
 \rangle 
 {{g}}' ( \delta ( u ) ) 
 \big] 
 \\ 
 & =  
 E\big[ 
 G \langle 
 ( \widetilde{\nabla} u )^k 
 u , 
 u 
 \rangle 
 {{g}} ( \delta ( u ) ) 
 \big] 
+ 
 E \big[ 
 \langle \widetilde{\nabla}^* u , 
 D ( 
 G 
 {{g}} ( \delta (u) ) 
 ( \widetilde{\nabla} u)^k u 
 ) 
 \big]
 \rangle 
 - E\big[ 
 G \langle 
 ( \widetilde{\nabla} u )^{k+1} 
 u  
 , 
 D 
 \delta ( u ) 
 \rangle 
 {{g}}' ( \delta ( u ) ) 
 \big] 
\\ 
 & =  
 E\big[ 
 G \langle 
 ( \widetilde{\nabla} u )^k 
 u , 
 u 
 \rangle 
 {{g}} ( \delta ( u ) ) 
 \big] 
 + 
 E \big[ 
 \langle (\widetilde{\nabla} u)^{k+1} u , 
 D G \rangle 
 {{g}} ( \delta (u) ) 
 \big] 
 + 
 E \big[ 
   G
 \langle \widetilde{\nabla}^* u , 
 D ( ( \widetilde{\nabla} u)^k u ) 
 \rangle 
   {{g}} ( \delta (u) ) 
 \big] 
\\ 
 & =  
 E\big[ 
 {{g}} ( \delta ( u ) ) 
 \Gamma_{k+2}^u G 
 \big],  
\end{align*} 
where we used \eqref{comp2} and \eqref{comp3}.
 Therefore, we have  
\begin{align*} 
 & 
 E[ 
 G 
 \delta ( u ) 
 {{g}} ( \delta ( u ) ) 
 ] 
 = 
 E[ 
   \langle u , D ( G {{g}} ( \delta ( u ) ) )
   \rangle 
 ] 
\\ 
 & =  
 E[ 
 G 
 \langle u , D \delta ( u ) \rangle 
 {{g}}' ( \delta ( u ) ) 
 ] 
 + 
 E[ 
 \langle u , D G \rangle 
 {{g}} ( \delta ( u ) ) 
 ] 
\\ 
 & =  
 E[ 
 \langle u , D G \rangle 
 {{g}} ( \delta ( u ) ) 
 ] 
 + 
 E\big[ 
 G  {{g}}^{(n+1)} ( \delta ( u ) ) 
\langle 
 ( \widetilde{\nabla} u )^n 
 u  
 , 
 D 
 \delta ( u ) 
 \rangle 
 \big] 
\\ 
 &  
 \quad + 
 \sum_{k=0}^{n-1} 
 \big( 
 E\big[ 
 G 
 {{g}}^{(k+1)} ( \delta ( u ) ) 
 \langle 
 ( \widetilde{\nabla} u )^k 
 u , 
 D 
 \delta ( u )
 \rangle 
 \big] 
 - 
 E\big[ 
 G  {{g}}^{(k+2)} ( \delta ( u ) ) 
\langle 
 ( \widetilde{\nabla} u )^{k+1} 
 u  
 , 
 D 
 \delta ( u ) 
 \rangle 
 \big] 
 \big) 
\\ 
 & =  
 E[ 
 \langle u , D G \rangle 
 {{g}} ( \delta ( u ) ) 
 ] 
 + 
 \sum_{k=1}^n 
 E\big[ 
 {{g}}^{(k)} ( \delta ( u ) ) 
 \Gamma_{k+1}^u G 
 \big] 
 + 
 E\big[ 
 G  {{g}}^{(n+1)} ( \delta ( u ) ) 
\langle 
 ( \widetilde{\nabla} u )^n 
 u  
 , 
 D 
 \delta ( u ) 
 \rangle 
 \big] 
\\ 
 & =  
 E[ 
 \langle u , D G \rangle 
 {{g}} ( \delta ( u ) ) 
 ] 
 + 
 \sum_{k=1}^n 
 E\big[ 
 {{g}}^{(k)} ( \delta ( u ) ) 
 \Gamma_{k+1}^u G 
 \big] 
 \\
 & \quad
 + 
 E\big[ 
 G  {{g}}^{(n+1)} ( \delta ( u ) ) 
\langle ( \widetilde{\nabla} u )^n u , u \rangle 
 \big] 
 + 
 E\big[ 
 G  {{g}}^{(n+1)} ( \delta ( u ) ) 
\langle 
 ( \widetilde{\nabla} u )^n 
 u  
 , 
 \delta ( \widetilde{\nabla}^* u ) 
 \rangle 
 \big], 
\end{align*} 
 and we conclude by Lemma~\ref{itlw2.1}.
\end{Proof} 
When $f\in W^{1,1}_0(B(R))$
is a deterministic function such that
$\Vert \nabla^{\real^d} f\Vert_\infty < \infty$, and $g\in {\cal C}^\infty_b ( \real )$, 
Proposition~\ref{djklsad1} shows that 
\begin{eqnarray*} 
  \lefteqn{
    E\left[ 
 \delta ( f ) 
 {{g}} ( \delta ( f ) ) 
 \right] 
  }
  \\
  & = & 
 \sum_{k=1}^{n+1}
 \frac{1}{k!}
 \int_{B(R)} f^{k+1}(x) \ \! \lambda ( dx ) 
 E[ 
 {{g}}^{(k)} ( \delta ( f ) ) 
 ] 
 + E\big[ 
 {{g}}^{(n+1)} ( \delta ( f ) ) 
 \langle 
 ( \widetilde{\nabla} f )^n 
 f  
 , 
 \delta ( \widetilde{\nabla}^* f ) 
 \rangle 
 \big] 
 \\
  & = & 
 \sum_{k=1}^{n+1}
 \frac{1}{k!}
 \kappa_{k+1}^f 
 E[ 
 {{g}}^{(k)} ( \delta ( f ) ) 
 ] 
 + E\big[ 
 {{g}}^{(n+1)} ( \delta ( f ) ) 
 \delta ( ( \widetilde{\nabla} f )^{n+1} f ) 
 \big], \qquad n \geq 0,
\end{eqnarray*} 
 with, by Proposition~\ref{p1} applied with $p=2$ and $r=0$, 
\begin{eqnarray*} 
  E\big[ 
 | 
 \delta (
  ( \widetilde{\nabla} f )^{n+1} 
 f  
 ) 
 | 
 \big] 
 & \leq & 
 \sqrt{E\big[ 
 | 
 \delta ( ( \widetilde{\nabla} f )^{n+1} 
 f  
 ) 
 |^2 
 \big]} 
\\
 & = & 
 \Vert 
 ( \widetilde{\nabla} f )^{n+1} f  
 \Vert_{L^2(B(R))}
\\
 & \leq & 
  ( K_d v_d R' )^{n+1} 
  \Vert f \Vert_{L^2 (B(R))}
  \Vert \widetilde{\nabla} f \Vert_{L^\infty (B(R);\real^d)}^{n+1}. 
\end{eqnarray*} 
 In addition, as $n$ tends to $+\infty$ we have 
\begin{eqnarray*} 
 E\left[ 
 \delta ( f ) 
 {{g}} ( \delta ( f ) ) 
 \right] 
 & = & 
 \sum_{k=1}^\infty 
  \frac{1}{k!}
 \int_{B(R)} 
f^{k+1}(x) 
 \ \! \lambda ( dx ) 
 E\left[ 
 {{g}}^{(k)} ( \delta ( f ) ) 
 \right] 
\\
\nonumber
& = & 
 \sum_{k=1}^\infty 
 \frac{1}{k!}
 \int_{B(R)} 
 f^{k+1}(x) \ \! \lambda ( dx ) 
 E\left[ 
 {{g}}^{(k)} ( \delta ( f ) ) 
 \right] 
 \\
  & = & E\left[
 \int_{B(R)} f(x) 
 \big( {{g}} ( \delta ( f ) + f(x) ) - {{g}} ( \delta ( f ) ) \big) 
 \lambda ( dx ) 
 \right]
\end{eqnarray*} 
 provided that the derivatives of $g$ decay fast enough, 
 which is a particular instance of
the standard integration by parts identity 
for finite difference operators on the Poisson space,
see e.g. Lemma~2.9 in \cite{utzet2} or Lemma~5 in \cite{bourguin}. 
\section{Stein approximation} 
\label{s4} 
Applying Proposition~\ref{djklsad1}
with $n=0$ and $G=1$ to the solution $g_x$ of the Stein equation 
\begin{equation} 
\nonumber 
 {\bf 1}_{(-\infty , x ]}(z) - \Phi ( z ) = g'_x ( z ) - z g_x (z), \qquad z\in\real,
\end{equation} 
 and $u \in \widetilde{\dee}^{1,1}_0$ a predictable random field 
 this gives the expansion 
\begin{eqnarray} 
\nonumber 
 P ( \delta ( u ) \leq x ) 
 - 
 \Phi ( x ) 
 & = & 
 E\big[ {{g}}'_x ( \delta ( u ) ) 
 \langle 
 u  
 , 
 u \rangle 
 - 
 \delta ( u ) 
 {{g}}_x ( \delta ( u ) ) 
 \big] 
 \\
 \nonumber
 & = & 
 E\left[ 
 ( 1 - \langle 
 u  
 , 
 u \rangle 
 ) {{g}}'_x ( \delta ( u ) ) 
\right] 
 + 
 E\big[ 
 \langle u , \delta ( \widetilde{\nabla} u ) \rangle 
 {{g}}'_x ( \delta ( u ) ) 
 \big],
\end{eqnarray} 
 around the Gaussian cumulative distribution function $\Phi (x)$,
 with $\Vert g_x \Vert_\infty \leq \sqrt{2\pi}/4$ 
 and $\Vert g'_x \Vert_\infty \leq 1$,
 $x\in \real$, by Lemma~2.2-$(v)$ of \cite{chenbk}. 
 The next result applies Proposition~\ref{djklsad1}
 with $n=1$ and $G=1$. 
\begin{prop} 
\label{djhlkdsa} 
For any random field $u \in \widetilde{\dee}^{1,\infty}_0$ 
 we have 
\begin{eqnarray} 
\nonumber
\lefteqn{
  d ( \delta ( u ) , {\cal N} ) 
}
\\
\nonumber
& \leq & 
 E\big[ 
 | 1 -  \langle 
 u 
 , 
 u 
 \rangle 
 -  \langle \widetilde{\nabla}^* u , D u \rangle 
 | 
 \big] 
+ 
   E\left[ 
   \left|
   \int_{B(R)} u_x^3 \ \! \lambda ( dx ) 
 + 
\left< 
 u , 
 D 
 \int_{B(R)} u^2_x \ \! \lambda ( dx ) 
 \right> 
 \right|
 \right] 
\\
\label{djklsadsa11} 
  & &
 + 
 2 E\big[ 
 | 
 \langle 
 ( \widetilde{\nabla} u ) 
 u  
 , 
 \delta ( \widetilde{\nabla}^* u ) 
 \rangle 
 | 
 \big] 
. 
\end{eqnarray} 
\end{prop} 
\begin{Proof} 
 For $n=1$ and $G=1$, Proposition~\ref{djklsad1} shows that 
\begin{eqnarray*} 
\nonumber 
 E [ \delta ( u ) 
 {{g}} ( \delta ( u ) ) 
 ] 
 & = & 
 E [ 
 {{g}}' ( \delta ( u ) ) 
 ( \langle u , u \rangle 
 + \langle \widetilde{\nabla}^* u , D u \rangle 
 ) 
 ] 
 \\
 & &
 + 
 \frac{1}{2} 
 E\left[ 
 g'' ( \delta ( u ) ) 
 \left(
 \int_{B(R)} u_x^3 \ \! \lambda ( dx ) 
 + 
 \left< 
 u , 
 D 
 \int_{B(R)} u^2_x \ \! \lambda ( dx ) 
 \right> 
 \right)
 \right] 
\\ 
 & & 
 + 
 E [ 
 {{g}}'' ( \delta ( u ) ) 
 \langle 
 ( \widetilde{\nabla} u ) 
 u  
 , 
 \delta ( \widetilde{\nabla} u ) 
 \rangle 
 ] 
. 
\end{eqnarray*} 
Let $h:\real\to [0,1]$
be a continuous function 
with bounded derivative. 
Using the solution $g_h\in\mathcal{C}_b^1(\real)$ 
 of the Stein equation
\begin{equation} 
\nonumber 
 h( z )-\mathrm{E}[h( {\cal N}  )] = g' ( z )- z g( z ), \qquad z\in\real,
\end{equation} 
 with the bounds
 $\Vert g'_h \Vert_\infty \leq \Vert h' \Vert_\infty$
 and
 $\Vert g''_h \Vert_\infty 
 \leq 2 \Vert h' \Vert_\infty$, 
 $x \in \real$, 
 cf. Lemma~1.2-$(v)$ of \cite{nourdinpeccati} and references therein,
 we have 
\begin{eqnarray*} 
\nonumber 
 E [ h ( \delta ( u ) ) ] 
 - 
 E [ h ( {\cal N} ) ] 
 & = & 
 E [ 
 \delta ( u ) 
 g_h ( \delta ( u ) ) 
 - 
 g_h' ( \delta ( u ) ) 
 ] 
\\ 
 & = & 
 E [ 
 g_h' ( \delta ( u ) ) 
 ( \langle u , u \rangle 
 +  \langle \widetilde{\nabla}^* u , D u \rangle 
 - 1 ) 
 ] 
 \\
  & & + 
 \frac{1}{2} 
 E\left[ 
 g'' ( \delta ( u ) ) 
 \left(
 \int_{B(R)} u_x^3 \ \! \lambda ( dx ) 
 + \left< 
 u , 
 D 
 \int_{B(R)} u^2_x \ \! \lambda ( dx ) 
 \right> 
 \right) \right] 
 \\
 & &
 + 
 2 E[ 
 g_h'' ( \delta ( u ) ) 
 \langle 
 ( \widetilde{\nabla} u ) 
 u  
 , 
 \delta ( \widetilde{\nabla}^* u ) 
 \rangle 
 ] 
, 
\end{eqnarray*} 
 hence 
\begin{eqnarray*} 
\nonumber 
  | E [ \delta ( u ) 
 h ( \delta ( u ) ) 
 ] 
 - 
 E [ 
 h ( {\cal N} ) 
 ] | 
& \leq & 
 \Vert h' \Vert_\infty 
 E\big[ 
   | 1 - \langle u , u \rangle 
   -  \langle \widetilde{\nabla}^* u , D u \rangle 
| 
 \big] 
 \\
   &  &
 + 
 \Vert h' \Vert_\infty 
 E\left[ 
   \left|
   \int_{B(R)} u_x^3 \ \! \lambda ( dx ) 
 + 
\left< 
 u , 
 D 
 \int_{B(R)} u^2_x \ \! \lambda ( dx ) 
 \right> 
 \right|
 \right] 
 \\
   &  &
   + 
 2 \Vert h' \Vert_\infty 
 E\big[ 
 | 
 \langle 
 ( \widetilde{\nabla} u ) 
 u  
 , 
 \delta ( \widetilde{\nabla}^* u ) 
 \rangle 
 | 
 \big] 
, 
\end{eqnarray*} 
which yields \eqref{djklsadsa11}.
\end{Proof} 
 As a consequence of Proposition~\ref{djhlkdsa} 
 and the It\^o isometry \eqref{isomrel} we have the following corollary.  
\begin{corollary} 
\label{djkd}
For $u \in \widetilde{\dee}^{1,\infty}_0$
 we have 
\begin{eqnarray} 
 \nonumber
 d ( \delta ( u ) , {\cal N} ) 
 & \leq & 
 | 1 - \Var [ \delta (u) ] | 
 + 
 \sqrt{ 
 \Var \big[ 
 \Vert u \Vert_{L^2(B(R))}^2 
 \big] 
 } 
 \\
 \nonumber 
 & &
 + 
   E\left[ 
   \left|
   \int_{B(R)} u_x^3 \ \! \lambda ( dx ) 
 + 
\left< 
 u , 
 D 
 \int_{B(R)} u^2_x \ \! \lambda ( dx ) 
 \right> 
 \right|
 \right] 
\\
\nonumber 
 & &
+ E[| \langle \widetilde{\nabla}^* u , D u \rangle 
  | ] 
+ 
 2 E\big[ 
 | 
 \langle 
 ( \widetilde{\nabla} u ) 
 u  
 , 
 \delta ( \widetilde{\nabla}^* u ) 
 \rangle 
 | 
 \big] 
. 
\end{eqnarray} 
\end{corollary} 
\begin{Proof} 
  By the It\^o isometry \eqref{isomrel} we have
$$ 
  \Var [ \delta (u) ] = E \left[ \left(
    \int_{B(R)} u_x ( \gamma ( dx ) - \lambda (dx) ) \right)^2 \right] 
 = 
 E [ 
 \langle 
 u 
 , 
 u 
 \rangle 
 ], 
 $$
  hence
\begin{eqnarray*} 
\nonumber
\lefteqn{
  E\big[ 
 | 1 -  \langle 
 u 
 , 
 u 
 \rangle 
-  \langle \widetilde{\nabla}^* u , D u \rangle 
 | 
 \big] 
}
\\
& \leq &
 E\left[ 
 | 1  -
 E [ \langle 
 u 
 , 
 u 
 \rangle] 
 | 
 \right] 
 +
   E\left[ 
 | \langle 
 u 
 , 
 u 
 \rangle 
 -
 E [ \langle 
 u 
 , 
 u 
 \rangle] 
 | 
 \right] 
   + E[| \langle \widetilde{\nabla}^* u , D u \rangle 
     | ] 
\\
 & = &
 | 1 - \Var [ \delta ( u ) ] | 
 +
 \sqrt{E [ ( \langle u , u \rangle 
     - E[\langle u , u \rangle 
     ] )^2 ]} 
 + E[| \langle \widetilde{\nabla}^* u , D u \rangle 
   | ] 
 \\
  & = & 
 | 1 - \Var [ \delta ( u ) ] | 
 +
  \sqrt{\Var \big[ \Vert u \Vert_{L^2(B(R))}^2 \big] } 
  + E[| \langle \widetilde{\nabla}^* u , D u \rangle 
    | ] 
. 
\end{eqnarray*} 
\end{Proof}
In particular, when $\Var [ \delta (u) ] = 1$,
Corollary~\ref{djkd} 
shows that 
\begin{eqnarray*} 
 d ( \delta ( u ) , {\cal N} ) 
 & \leq & 
 \sqrt{ 
 \Var \big[ 
 \Vert u \Vert_{L^2(B(R))}^2 
 \big] }
 + 
   E\left[ 
   \left|
   \int_{B(R)} u_x^3 \ \! \lambda ( dx ) 
 + 
\left< 
 u , 
 D \int_{B(R)} u^2_x \ \! \lambda ( dx ) 
 \right> 
 \right|
 \right] 
   \\
   & &
   + E[| \langle \widetilde{\nabla}^* u , D u \rangle 
     | ] 
 + 
 2 E\big[ 
 | 
 \langle 
 ( \widetilde{\nabla} u ) 
 u  
 , 
 \delta ( \widetilde{\nabla}^* u ) 
 \rangle 
 | 
 \big] 
.  
\end{eqnarray*} 
 When $f\in W^{1,\infty}_0(B(R))$
is a deterministic function we have 
$$
 \Var [ \delta (f) ]
=
\E \left[
  \left(
  \int_{B(R)} f(x) ( \gamma ( dx ) - \lambda (dx) )  
  \right)^2
  \right]
=
\int_{B(R)} f^2(x) \ \! \lambda ( dx ), 
$$
 and Corollary~\ref{djhlkdsa} shows that 
$$
 d ( \delta ( f ) , {\cal N} ) 
 \leq 
   \left| 1 - \int_{B(R)} f^2(x) \ \! \lambda (dx) 
 \right| 
  + 
  \left|
  \int_{B(R)} f^3(x) \ \! \lambda ( dx ) 
 \right|
 +
2  E\big[ 
 | 
 \delta ( ( \widetilde{\nabla} f )^2 f ) 
 | 
 \big] 
. 
$$ 
 Given the bound 
\begin{eqnarray*} 
  E\big[ 
 | 
 \delta (
  ( \widetilde{\nabla} f )^2 
 f  
 ) 
 | 
 \big] 
 & \leq & 
 \sqrt{E\big[ 
 | 
 \delta (
  ( \widetilde{\nabla} f )^2 
 f  
 ) 
 |^2 
 \big]} 
\\
 & = & 
 \Vert 
   ( \widetilde{\nabla} f )^2 
 f  
 \Vert_{L^2(B(R))}
\\
 & \leq & 
  ( K_d v_d R' )^2 
  \Vert f \Vert_{L^2 (B(R))}
  \Vert \nabla^{\real^d} f \Vert_{L^\infty (B(R);\real^d)}^2  
\end{eqnarray*} 
obtained from Proposition~\ref{p1} with $p=2$ and $r=0$,
$f\in W^{1,\infty}_0(B(R))$, we have the following corollary.
\begin{corollary} 
 For $f\in W^{1,\infty}_0(B(R))$ we have 
\begin{eqnarray} 
\nonumber 
 d \left( \int_{B(R)} f(x) ( \gamma ( dx ) - \lambda (dx) ), {\cal N} \right) 
 & \leq & 
   \big| 1 - \Vert f\Vert_{L^2(B(R))}^2 
 \big| 
+ \left| \int_{B(R)} f^3(x) \ \! \lambda ( dx ) \right| 
\\
\nonumber 
& &
+ 2 ( K_d v_d R' )^2 
  \Vert f \Vert_{L^2 (B(R))}
  \Vert \nabla^{\real^d} f \Vert_{L^\infty (B(R);\real^d)}^2. 
\end{eqnarray} 
\end{corollary} 
In particular, if $\Vert f\Vert_{L^2(B(R))}=1$ we find 
$$ 
 d \left( \int_{B(R)} f(x) ( \gamma ( dx ) - \lambda (dx) ), {\cal N} \right) 
 \leq
 \left| \int_{B(R)} f^3(x) \ \! \lambda ( dx ) \right| 
+ 2 ( K_d v_d R' )^2 
\Vert \nabla^{\real^d} f \Vert_{L^\infty (B(R);\real^d)}^2.
$$
As an example, consider $f_k$ given on
$B(Rk^{1/d})$ by 
$$
f_k (x) := \frac{1}{C \sqrt{k} } g\left(\frac{|x|_{\real^d}}{k^{1/d}}\right)
, 
$$
where $g\in {\cal C}^1_0 ([0,R])$ 
and 
$$
C^2 := v_d \int_0^R g^2(r) r^{d-1} dr,
$$
so that
$f_k \in L^2(B(R k^{1/d}))$ with 
$$
 \Vert f \Vert_{L^2(B(R k^{1/d}))}^2 
=
\frac{v_d}{C^2 k}
\int_0^{ R k^{1/d}} g^2 \left( \frac{r}{k^{1/d}} \right) r^{d-1} dr
=
\frac{v_d}{C^2} \int_0^R g^2 (r) r^{d-1} dr
=
1,
$$
 and
$$
\int_{B(Rk^{1/d})} f_k^3 (x) dx
=
\frac{1}{C^3 k^{3/2}}
\int_0^{Rk^{1/d}} g^3 (r k^{-1/d}) r^{d-1} dr
=
\frac{1}{C^3\sqrt{k}}
\int_0^R g^3 (r) r^{d-1} dr, 
$$
 $k \geq 1$.
We have
$$
\Vert \nabla^{\real^d} f_k \Vert_{L^\infty (B(R);\real^d)}^2 \leq
\frac{\Vert g' \Vert_\infty^2 d}{C^2 k^{1+2/d}}, 
$$
 hence 
\begin{eqnarray*} 
 d \left( \int_{B(R)} f_k(x) ( \gamma ( dx ) - \lambda (dx) ), {\cal N} \right) 
 & \leq & 
 \left| \int_{B(R)} f^3_k(x) \ \! \lambda ( dx ) \right| 
 + \frac{ 2 ( K_d v_d k^{1/d} R' )^2 d }{ k^{1+2/d} C^2}\Vert g' \Vert_\infty^2
\\
 & \leq & 
\frac{v_d}{C^3\sqrt{k}}
\left| \int_0^R g^3 (r) r^{d-1} dr \right| 
+ \frac{ 2 ( K_d v_d R' )^2 d }{ k C^2}\Vert g' \Vert_\infty^2.
\end{eqnarray*} 
In particular, if $g$ satisfies the condition
$$
\int_0^R g^3 (r) r^{d-1} dr =0, 
$$
 then we find the $O(1/k)$ convergence rate
$$ 
 d \left( \int_{B(R)} f_k(x) ( \gamma ( dx ) - \lambda (dx) ), {\cal N} \right) 
 \leq \frac{ 2 ( K_d v_d R' )^2 d}{ k C^2},
 \qquad k \geq 1. 
$$ 
 \section{Appendix}
 \begin{Proofy} \hskip-0.32cm {\em of Proposition~\ref{prop1}.}
       
 As a consequence of \eqref{gxy4} and \eqref{cont} we have
\begin{eqnarray} 
\nonumber
\lefteqn{
 f_n \left( x_1,\ldots , x_{i-1}, y , x_{i+1}, \ldots x_n \right)
 -
 f_{n-1} \left( x_1,\ldots , x_{i-1},x_{i+1}, \ldots, x_n \right)
}
\\
\nonumber
& = &
 f_n \left( x_1,\ldots , x_{i-1}, y , x_{i+1}, \ldots x_n \right)
 -
 f_{n-1} \left( x_1,\ldots , x_{i-1},x_{i+1}, \ldots, x_n \right)
 \int_{B(R')\setminus B(R)}
 \eta(x)
 \ \! \lambda (dx)
 \\
\nonumber
  & = & 
 f_n \left( x_1,\ldots , x_{i-1}, y , x_{i+1}, \ldots x_n \right)
 -
 \int_{B(R')\setminus B(R)}
 \eta(x)
 f_n \left( x_1,\ldots , x_{i-1},x,x_{i+1}, \ldots, x_n \right)
 \ \! \lambda (dx)
 \\
\nonumber 
  & = & 
\int_{B(R')} \langle G(x_i,y)
,
\nabla^{\real^d}_{x_i} f_n \left( x_1,\ldots , x_n \right)
\rangle_{\real^d} 
\lambda ( dx_i )
 \\
\label{cont2} 
  & = & 
\int_{B(R)} \langle G(x_i,y)
,
\nabla^{\real^d}_{x_i} f_n \left( x_1,\ldots , x_n \right)
\rangle_{\real^d} 
\lambda ( dx_i )
,
\end{eqnarray} 
$x_1,\ldots , x_{i-1},y,x_{i+1}, \ldots , x_n\in B(R')$.
 Recall that for all $F\in {\mathcal S}$ of the form \eqref{f} we have
$$ 
 E [F] = e^{-B(R)}
 \sum_{n=0}^\infty \frac{1}{n!} 
 \int_{B(R)} \cdots \int_{B(R)} f_n(x_1,\ldots ,x_n) \ \! \lambda ( dx_1 ) \cdots \lambda (dx_n).
$$
 Hence, using \eqref{cont2},
 for $g\in {\cal C}^1_0 (B(R))$ and $F$ of the form \eqref{f} 
 we have 
  \begin{align} 
  \nonumber
  & 
  E \left[ \int_{B(R)} g(y) D_y F \ \! \lambda ( dy ) \right]
 \\
  & = 
 E\left[
   \sum_{n=1}^\infty 
   {\bf 1}_{\{ \gamma ( B(R))=n\} }
 \sum_{i=1}^n 
\int_{B(R)}
 g(y) \langle \G_\eta (X_i,y) ,
 \nabla^{\real^d}_{X_i} f\left( X_1,\ldots , X_n \right)
 \rangle \lambda ( dy)
 \right]
    \\
  \nonumber
  & = 
  e^{-B(R)}
 \sum_{n=1}^\infty \frac{1}{n!} 
 \int_{B(R)} \cdots \int_{B(R)}
 \sum_{i=1}^n 
\int_{B(R)}
 g(y) \langle \G_\eta (x_i,y) ,
 \nabla^{\real^d}_{x_i} f_n(x_1,\ldots ,x_n)
 \rangle \lambda ( dy)
   \lambda ( dx_1 ) \cdots \lambda (dx_n)
    \\
  \nonumber
  & = 
  e^{-B(R)}
  \sum_{n=1}^\infty
  \frac{1}{n!} 
  \int_{B(R)} \cdots
  \\
  \nonumber
  & \quad 
\cdots \int_{B(R)}
 \sum_{i=1}^n 
\int_{B(R)}
g(y)
 f_n(x_1,\ldots ,x_{i-1},y,x_{i+1},\ldots , x_n)
   \ \! \lambda ( dx_1 ) \cdots  \lambda ( dy) \cdots \lambda (dx_n)
   \\
   \nonumber 
        & 
 \quad - 
  e^{-B(R)}
 \sum_{n=1}^\infty \frac{1}{n!} 
 \int_{B(R)} \cdots \int_{B(R)}
 \sum_{i=1}^n 
\int_{B(R)}
g(y)
 \ \! \lambda ( dy)
 f_{n-1}(x_1,\ldots ,x_{n-1})
   \ \! \lambda ( dx_1 ) \cdots \lambda (dx_{n-1})
    \\
  \nonumber
  & = 
  e^{-B(R)}
  \sum_{n=1}^\infty
  \frac{1}{n!} 
  \int_{B(R)} \cdots \int_{B(R)}  
  \left(
  \sum_{i=1}^n 
g(x_i)
- \int_{B(R)}
g(y)
 \ \! \lambda ( dy)
 \right)
 f_n (x_1,\ldots ,x_n )
 \ \! \lambda ( dx_1 ) \cdots \lambda (dx_n)
    \\
\nonumber 
  & =   
    E\left[
      F
      \left(
      \int_{B(R)} g(x) ( \gamma ( dx ) - \lambda (dx) ) 
      \right)
      \right]. 
\end{align} 
Next, for $u$ of the form \eqref{u}, we check by a standard argument that
\begin{eqnarray*} 
  \nonumber
E[    \langle u , DF\rangle_{L^2(B(R))} ] 
  & = & 
\sum_{i=1}^n
E[ G_i \langle g_i , DF\rangle_{L^2(B(R))} ] 
\\
& = & 
\sum_{i=1}^n
\left(
E[ \langle g_i , D( F G_i ) \rangle_{L^2(B(R))} - F \langle g_i , DG_i \rangle_{L^2(B(R))} ] 
\right)
\\
& = & 
    E\left[
      F
\sum_{i=1}^n
\Big(
G_i 
     \int_{B(R)} g_i(x) ( \gamma ( dx ) - \lambda (dx) ) 
- 
\langle g_i , DG_i \rangle_{L^2(B(R))}
\Big)
    \right] 
\\
& = & 
    E [ 
      F \delta ( u )].
    \end{eqnarray*} 
\end{Proofy} 
\begin{Proofy} \hskip-0.32cm {\em of Proposition~\ref{prop2}.}  
 Taking $u \in {\cal P}_0$ a predictable random field
of the form \eqref{uf} we note that by \eqref{ndfgd} and the compatibility condition
\eqref{compat} we have 
$$
g_i(y)D_yF_i=0, \qquad y\in B(R), \quad i=1,\ldots , n,
$$
hence by 
\eqref{dltu} we have 
\begin{eqnarray} 
  \label{dltu2} 
\delta ( u )
& = & \delta \left( \sum_{i=1}^n F_i g_i \right)
 = \sum_{i=1}^n F_i \delta ( g_i )
\\
\nonumber
& = & \sum_{i=1}^n F_i \int_{B(R)} g_i(x) ( \gamma ( dx ) - \lambda ( dx )) 
\\
\nonumber
& = & \int_{B(R)} u_x ( \gamma ( dx ) - \lambda ( dx )), 
\end{eqnarray} 
 showing that $\delta (u)$ coincides with the
 Poisson stochastic integral of $(u_x)_{x\in B(R)}$.
 Regarding the isometry relation \eqref{isomrel}, we have 
\begin{eqnarray*} 
      E [ \delta (u)^2 ]
 & = & E\left[
  \left(
  \sum_{i=1}^n
  F_i \int_{B(R)} g_i(x) ( \gamma ( dx ) - \lambda ( dx ) )
  \right)^2
  \right]
  \\
& = & 
E\left[
  \sum_{i,j=1}^n
  F_i F_j \int_{B(R)} g_i(x) ( \gamma ( dx ) - \lambda ( dx ) )
  \int_{B(R)} g_j(x) ( \gamma ( dx ) - \lambda ( dx ) )
  \right]
\\
& = & 
2 E\left[
  \sum_{1\leq i<j\leq n}
  F_i \int_{B(R)} g_i(x) ( \gamma ( dx ) - \lambda ( dx ) )
  F_j \int_{B(R)} g_j(x) ( \gamma ( dx ) - \lambda ( dx ) )
  \right]
\\
 & & +
E\left[
  \sum_{i=1}^n
  F_i^2 \left( \int_{B(R)} g_i(x) ( \gamma ( dx ) - \lambda ( dx ) )
  \right)^2 
  \right]
\\
& = & 
E\left[
  \sum_{i=1}^n
  F_i^2 \int_{B(R)} g_i^2 (x) \lambda ( dx ) 
  \right]
\\
& = & 
E\left[
\int_{B(R)} u^2 (x) \ \! \lambda ( dx ) 
  \right], 
\end{eqnarray*} 
which shows that \eqref{rel}
extends to the closure of ${\cal P}_0$ in $L^2(\Omega \times B(R))$
by density and a Cauchy sequence argument. 
\end{Proofy}

\footnotesize 
 
\def\cprime{$'$} \def\polhk#1{\setbox0=\hbox{#1}{\ooalign{\hidewidth
  \lower1.5ex\hbox{`}\hidewidth\crcr\unhbox0}}}
  \def\polhk#1{\setbox0=\hbox{#1}{\ooalign{\hidewidth
  \lower1.5ex\hbox{`}\hidewidth\crcr\unhbox0}}} \def\cprime{$'$}

\end{document}